\documentclass[10pt,reqno]{article}
\usepackage{amsmath,amssymb,amsthm}
\usepackage{esint}
\usepackage{bm}
\usepackage{url}
\usepackage{subfigure}
\usepackage{graphicx,color}
\usepackage{algorithm}
\usepackage{algpseudocode}
\usepackage{algorithmicx}
\usepackage{hyperref}
\usepackage{float}
\usepackage{booktabs,multirow}
\usepackage{hhline}
\usepackage{setspace}

\usepackage{caption}
\usepackage{hyphsubst}
\usepackage{caption}

\usepackage{float}
\usepackage{filecontents}
\usepackage{hyperref}

\hypersetup{colorlinks,citecolor=blue,linkcolor=blue, urlcolor=blue}
\usepackage{pdflscape}
\usepackage{breakurl}

\newcommand{\RR}{\mathbb{R}}

\DeclareMathAlphabet{\itbf}{OML}{cmm}{b}{it}

\newtheorem{thm}{Theorem}[section]

\newtheorem{rem}[thm]{Remark}

\numberwithin{equation}{section}

\newcommand{\email}[1]{\protect\href{mailto:#1}{#1}}

\newcommand{\pathfigures}{Figures/}
\graphicspath{{\pathfigures}}

\begin{document}

\title{Comments on \textquotedblleft Design of momentum fractional LMS for Hammerstein nonlinear system identification with application to electrically stimulated muscle model\textquotedblright
}

\author{
Abdul Wahab\footnotemark[1]\, \footnotemark[2]
\and
Shujaat Khan\footnotemark[3]
\and 
Farrukh Zeeshan Khan\footnotemark[4]
}
\maketitle
\renewcommand{\thefootnote}{\fnsymbol{footnote}}
\footnotetext[1]{Corresponding Author. E-mail address:  \email{abdul.wahab@nu.edu.kz}.}
\footnotetext[2]{Department of Mathematics, School of Sciences and Humanities, Nazarbayev University, 53 Kabanbay Batyr Avenue, Nur-Sultan 010000, Kazakhstan (\email{abdul.wahab@nu.edu.kz)}.}
\footnotetext[3]{Bio-Imaging, Signal Processing, and Learning Lab., Department of Bio and Brain Engineering, Korea Advanced Institute of Science and Technology, 291 Daehak-ro, Yuseong-gu, 34141, Daejeon, South Korea (\email{shujaat@kaist.ac.kr}).}
\footnotetext[4]{Department of Computer Science, University of Engineering and Technology Taxila, 47080, Taxila, Pakistan (\email{farrukh.zeeshan@uettaxila.edu.pk)}.}
\renewcommand{\thefootnote}{\arabic{footnote}}

\begin{abstract}
The purpose of this article is to discuss some aspects of the convergence analysis performed in the paper [Design of momentum fractional LMS for Hammerstein nonlinear system identification with application to electrically stimulated muscle model, Eur. Phys. J. Plus (2019) \textbf{134}: 407]. It is highlighted that the way the authors prove convergence suffers a lack of correct and valid mathematical justifications. 
\end{abstract}



\section{Introduction}\label{S:Intro}

In \cite{paper}, the so-called \textit{momentum fractional least-mean-square (M-FLMS)} algorithm is utilized for estimation of parameters of a nonlinear system based on Hammerstein autoregressive exogenous input (H-ARX) structure. The exploited M-FLMS algorithm was first introduced in \cite{MFLMS} whereas the Hammerstein nonlinear control autoregressive systems were discussed earlier using another fractional variant of the LMS algorithm in \cite{ref37}. The main contribution of \cite{paper} is the mathematical and computational performance analysis of the M-FLMS algorithm presented in \cite{MFLMS} for a stimulated muscle modeling and identification system problem. The most important part of the study is the mathematical analysis which is presented in \cite[Sec. 4]{paper}. However, there are mathematically invalid justifications that are detrimental to the correctness of the entire framework. The main objective of this note is to list those mistakes and substantiate that the entire analysis is mathematically invalid.

\begin{rem}
The symbols, notations, and equation numbers used in this comment are consistent with \cite{paper}.
\end{rem}

\section{Mathematical Mistakes}\label{S:Errors}
Let us first give some general remarks regarding the overall design of the problem undertaken in \cite{paper}.

\subsection{General Remarks}

\begin{enumerate}
\item Basic matrix algebra suggests that \cite[Eq. 8]{paper} is invalid. Indeed, each $\Psi_i\in\RR^m$ is being multiplied with $q_i\mathbf{c}\in\RR^l$ (dimension mismatch) and is therefore, invalid and different from \cite[Eq. 7]{paper}.  In fact, the definition of $\Psi_i\in\RR^m$ in \cite[Eq. 10]{paper} is faulty. For \cite[Eq. 8]{paper} to be valid and equivalent to \cite[Eq. 7]{paper}, $\Psi_i$ should have been defined as 
\begin{align}
\tag{10*}\label{eq10*}
\Psi_i(t) = \Big[f_1[r(t-i)], f_2[r(t-i)], \cdots, f_l[r(t-i)]\Big]^T\in\RR^l, \quad 1\leq i\leq m.
\end{align}

\item \cite[Eq. 17]{paper} is derived from \cite[Eq. 11]{paper} by applying rule \cite[Eq. 16]{paper} with $a=0$ and subsequently using ordinary chain rule.  It is a well-known fact that the chain rule for fractional derivatives is different and is much more complicated than that for ordinary derivatives, see e.g., \cite{Tarasov, Kilbas}. The expression \cite[Eq. 17]{paper} is mathematically unjustified. The fractional derivative \cite[Eq. 17]{paper} of the mean square error (MSE) \cite[Eq. 11]{paper} has already been discussed in detail in \cite{NODY} and the correct form is also presented therein. Nevertheless, the resulting expression \cite[Eq. 18]{paper} and the weight update rules \cite[Eq. 19-21]{paper} can be viewed as approximations and therefore, we don't make any claim.

\item \cite[Eq. 18]{paper} is a vector equation containing fractional power $\mathbf{\hat\Omega}^{1-v}(t)$ of a vector $\mathbf{\hat\Omega}(t)$. The equation  is taken from \cite{ref37} wherein the \textit{fractional power of a vector} is defined component-wise. Whenever, there is a negative element in the vector $\mathbf{\hat\Omega}$, $\mathbf{\hat\Omega}^{1-v}(t)$ will render a complex output and impede the algorithm to converge to a real sought value. See \cite{NODY, WS, Bershad} for detailed discussions on this.  In case $\mathbf{\hat\Omega}^{1-v}(t)$ is not defined item-wise, its sense needs to be specified since the fractional powers of vectors are mathematically undefined.

\item To avoid appearance of the complex outputs, a modulus is introduced in \cite[Eq. 19-21]{paper}, i.e.,  $\mathbf{\hat\Omega}^{1-v}(t)$ is replaced by $|\mathbf{\hat\Omega}|^{1-v}(t)$ (which makes it a scalar). However, an itemized vector multiplication operator $\otimes$ is introduced in the equations that does not make sense because now there is only one vector in each of the products $\mathbf{\Psi}(t)\otimes |\mathbf{\hat\Omega}|^{1-v}(t)$ and $\mathbf{\Psi}(t)\otimes(1+ |\mathbf{\hat\Omega}|^{1-v}(t))$ in \cite[Eq. 19-21]{paper}. 
\end{enumerate}

It is worthwhile mentioning that all the aforementioned errors can be taken as approximations in algorithmic design, yet, they can be nasty in convergence analysis. Therefore, we do not make any claim here and, instead, use this information to establish our claims regarding convergence analysis.

\subsection{Issues in Convergence Analysis}

\begin{enumerate}
\item The convergence analysis contained in \cite[Sect. 4.2]{paper} is based on \cite[Eq. 23]{paper} given by
\begin{align}
\Delta\mathbf{\Omega}(t+1)=&\Delta\mathbf{\Omega}(t)+\beta[\Delta\mathbf{\Omega}(t)-\Delta\mathbf{\Omega}(t-1)]
\nonumber
\\
&+\eta\mathbf{\Psi}(t)\Big\{s(t)-\mathbf{\Psi}(t)^T[\mathbf{\Omega}_{\rm opt}+\Delta\mathbf{\Omega}(t)]\Big\}\Big\{1+[\mathbf{\Omega}_{\rm opt}+\Delta\mathbf{\Omega}(t)]^{1-\upsilon}\Big\}.
\tag{23} \label{e23}
\end{align}
Remark that, the update rule for the M-FLMS is given by \cite[Eq. 21]{paper},
\begin{align}
\hat{\mathbf{\Omega}}(t+1)=\hat{\mathbf{\Omega}}(t)+\beta\left[\hat{\mathbf{\Omega}}(t)-\hat{\mathbf{\Omega}}(t-1)\right]
+\eta e(t)\mathbf{\Psi}(t)\otimes\Big\{1+\left|\hat{\mathbf{\Omega}}\right|^{1-\upsilon}(t)]\Big\},
\tag{21} \label{e21}
\end{align}
and contains $|\mathbf{\hat{\Omega}}|^{1-v}$, whereas the magnitude disappears in \cite[Eq. 23]{paper}. In fact, 
$|\mathbf{\hat{\Omega}}|^{1-v}$ is a scalar and $\mathbf{\hat{\Omega}}^{1-v}$ is a vector (with component-wise power). Therefore, \eqref{e21} does not provide \eqref{e23}. Moreover, vector $\left[\mathbf{\Omega}_{\rm opt}+\Delta\mathbf{\Omega(t)}\right]^{1-v}$ is added to a scalar $1$ in the last term of \eqref{e23} which is mathematically absurd. 

\item In equation \cite[Eq. 24]{paper}, vector $\mathbf{\Psi}(t)$ is multiplied with $\left[\mathbf{\Omega}_{\rm opt}+\Delta\mathbf{\Omega(t)}\right]^{1-v}$. Recall that, for all vectors $\textbf{u},\mathbf{v}\in\mathbb{R}^N$ for any positive integer $N$, the product $\mathbf{u}\mathbf{v}$ is not defined in the classical sense, whereas $\mathbf{u}\mathbf{v}^T$   gives us a dyad (or simply a matrix). Therefore, it is not clear how the product $\eta s(t)\mathbf{\Psi}(t)[\mathbf{\Omega}_{\rm opt}+\Delta\mathbf{\Omega}(t)]^{1-\upsilon}$ should be interpreted. Note that the resultant is then added to vectors $\Delta\mathbf{\Omega}(t)$,  and $\eta\mathbf{\Psi}(t)s(t)$ etc. in \cite[Eq. 24]{paper}, which does not make sense. 

\item A \textit{binomial} expansion given by (\cite[Eq. 25]{paper}),
\begin{align}
[\mathbf{\Omega}_{\rm opt}+\Delta\mathbf{\Omega}(t)]^j=\sum_{k=0}^\infty \begin{pmatrix}
j\\k
\end{pmatrix}
\left(\mathbf{\Omega}_{\rm opt}^k\right)^T\Delta\mathbf{\Omega}(t)^{j-k}, \qquad |\mathbf{\Omega}_{\rm opt}|>|\Delta\mathbf{\Omega}|,
\tag{25}
\end{align}
is used to derive \cite[Eq.27]{paper},
\begin{align}
\Delta\mathbf{\Omega}(t+1)=&\Delta\mathbf{\Omega}(t)+\beta[\Delta\mathbf{\Omega}(t)-\Delta\mathbf{\Omega}(t-1)]
\nonumber
\\
&+\eta\mathbf{\Psi}(t)s(t)-\eta\mathbf{\Psi}(t)\mathbf{\Psi}(t)^T\mathbf{\Omega}_{\rm opt}
-\eta\mathbf{\Psi}(t)\mathbf{\Psi}(t)^T\Delta\mathbf{\Omega}(t)
\nonumber
\\
&+\eta\mathbf{\Psi}(t)s(t)\Delta\mathbf{\Omega}(t)^{1-\upsilon}+\eta\mathbf{\Psi}(t)s(t)\sum_{k=1}^\infty \begin{pmatrix}
1-\upsilon\\k
\end{pmatrix}
\left(\mathbf{\Omega}_{\rm opt}^k\right)^T\Delta\mathbf{\Omega}(t)^{1-\upsilon-k}
\nonumber
\\
&-\eta\mathbf{\Psi}(t)\mathbf{\Psi}(t)^T\Delta\mathbf{\Omega}(t)^{2-\upsilon}
\nonumber
-\eta\mathbf{\Psi}(t)\mathbf{\Psi}(t)^T\begin{pmatrix}
2-\upsilon\\1
\end{pmatrix}
\left(\mathbf{\Omega}_{\rm opt}^{1}\right)^T\Delta\mathbf{\Omega}(t)^{1-\upsilon}
\\
&-\eta\mathbf{\Psi}(t)\mathbf{\Psi}(t)^T\sum_{k=1}^\infty \begin{pmatrix}
1-\upsilon\\k
\end{pmatrix}
\left(\mathbf{\Omega}_{\rm opt}^{k+1}\right)^T\Delta\mathbf{\Omega}(t)^{1-\upsilon-k},
\tag{27}\label{e27}
\end{align}
which is inappropriate for vectors. In fact, such a \textit{binomial} expansion is unwarranted in Mathematics. 
\begin{itemize} 
\item[(i)]
The basic reason is  that the multiplication of vectors, if possible, is not commutative. 

\item[(ii)] It is unclear how to  interpret  $[\mathbf{\Omega}^k_{\rm opt}]^T\Delta\mathbf{\Omega}(t)^{j-k}$? If the powers are defined component-wise then this product and the right-hand side of \cite[Eq. 25]{paper} is scalar, whereas, the left-hand side is a vector.  So, equality cannot hold.

\item[(iii)] How would one interpret $\mathbf{\Omega}^0_{\rm opt}$? A careful inspection of \cite[Eq. 27]{paper} demystifies that $\mathbf{\Omega}^0_{\rm opt}$ is taken to be $1$; see the terms $\eta\mathbf{\Psi}(t)s(t)\Delta\mathbf{\Omega}(t)^{1-\upsilon}$ and $-\eta\mathbf{\Psi}(t)\mathbf{\Psi}(t)^T\Delta\mathbf{\Omega}(t)^{2-\upsilon}$ (the first terms on the third and fourth lines of \eqref{e27}, respectively). However, with powers of vectors defined component-wise, $\mathbf{\Omega}^0_{\rm opt}$ should be a vector with all entries $1$.

\item[(iv)] The term $\eta\mathbf{\Psi}(t)s(t)\Delta\mathbf{\Omega}(t)^{1-\upsilon}$ (which has no clear mathematical sense) is incompatible with the rest of additive (vector) terms in \eqref{e27}. 

\end{itemize}

\item Given the above, \eqref{e27} is mathematically invalid and absurd. The rest of the convergence analysis is based on this equation. Moreover, \cite[Eq. 28]{paper} has similar issues.

\item It is claimed that $\mathbf{p}-\mathbf{R}\mathbf{\Omega}_{\rm opt}=0$ after \cite[Eq. 28]{paper}, where $\mathbf{\Omega}_{\rm opt}$ is the optimal Wiener solution. However, there is no guarantee that the fractional optimal solution $\mathbf{\hat{\Omega}}_{\rm opt}$, furnished by the M-FLMS algorithm, will converge to the Wiener solution or even if it will converge at all. We refer the reader to \cite{Minima} for a discussion on fractional extreme values and issues with fractional optimization frameworks in the family of M-FLMS.  

\item A function $\mathbf{F}(\Delta\mathbf{\Omega}(t), v)$ is defined in \cite[Eq. 30]{paper}. It is not clear whether it is a scalar field, vector field, or a matrix-valued function. The existence of such a function is questionable. We have the following remarks.

\begin{itemize}

\item[(i)] In the last term of \cite[Eq. 30]{paper}, i.e., $\eta E[\Delta\mathbf{\Omega}(t)]\mathbf{F}[\Delta\mathbf{\Omega}(t), v]$, the function $\mathbf{F}$ is multiplied with a column vector $\Delta\mathbf{\Omega}(t)$ from the left, which is only possible if $\mathbf{F}$ is a scalar. On contrary, it is being added to a matrix $\mathbf{R}$ in \cite[Eq. 32]{paper} which is valid only if it is a matrix. So, $\mathbf{F}$ needs to be a scalar and a matrix at the same time. 

\item[(ii)] From \cite[Eq. 32]{paper} to \cite[Eq. 37]{paper}, $\mathbf{F}$ is treated as a matrix of the same dimension as $\mathbf{R}$, this is apparent from the context, otherwise, all these equations are invalid. However, in \cite[Eq. 38]{paper}, it is once again subtracted from a scalar $\lambda_i$. Moreover, it is divided in \cite[Eq. 39]{paper}. However,division for matrices or vectors is not defined in Mathematics. So, either it is a scalar and \cite[Eqs. 32-37]{paper} are invalid or it is a matrix and \cite[Eqs 38-39]{paper} are invalid. Hence, for the convergence analysis to make some sense $\mathbf{F}$ can neither be a scalar nor a matrix. 
\end{itemize}
Based on the aforementioned observations, one can simply conclude that such a function $\mathbf{F}$ does not even exist.

\item Because of the above mathematical jargon, it is unclear how much the numerical values and the simulations provided in \cite[Sect. 5]{paper} are accurate. 

\end{enumerate}

Based on these observations, it is clear that the entire convergence analysis suffers a lack of correct and valid mathematical justifications.  

\section{Conclusions} 
In this note, we have provided details of mathematical mistakes in the convergence analysis presented in \cite{paper}. Based on our observations, it is concluded that the convergence analysis in \cite{paper} is invalid and mathematically wrong.

\section*{Conflict of Interest} 
The authors declare that they have no conflict of interest.

\bibliographystyle{plain}

\end{document}